\documentclass[11pt,reqno]{amsart}
\usepackage{mathrsfs}
\usepackage{amsgen}
\usepackage{amscd}
\usepackage{amsmath}
\usepackage{latexsym}
\usepackage{amsfonts}
\usepackage{amssymb}
\usepackage{amsthm}
\usepackage{graphicx}
\usepackage{color}
\usepackage{enumerate}
\usepackage{amssymb,amsmath,graphicx,amsfonts,euscript}
\usepackage{color}
\usepackage{empheq}
\usepackage{cases}
\usepackage{cite}

\usepackage[colorlinks=true,backref]{hyperref}
\hypersetup{urlcolor=blue, citecolor=green, linkcolor=blue}

\usepackage[margin=3cm, a4paper]{geometry}

\def\H{{\cal H}}

\def\N{\mathbb{N}}
\def\R{\mathbb{R}}
\def\Z{\mathbb{Z}}

\def\H2{H^2(\R^N)}
\def\L2{L^2(\R^N)}
\def\to{\rightarrow}

 \def\norm#1{\|#1\|}
\def\normb#1{\bigg\|#1\bigg\|} 

\def\H{{\cal H}}

\def\H1{H^1(\R)}

\newcommand{\al}{\alpha} 
 \newcommand{\de}{\delta}

 \newcommand{\Om}{\Omega}
\newcommand{\Ga}{\Gamma}

\newcommand{\re}{\mathop{\mathrm{Re}}}

 \newcommand{\Del}[1]{}

\numberwithin{equation}{section}

\newtheorem{thm}{Theorem}[section]

\newtheorem{lem}[thm]{Lemma}

\newtheorem{definition}[thm]{Definition}
\theoremstyle{remark}

\newtheorem*{exam*}{Examples}

\newcommand{\EQ}[1]{\begin{align*}\begin{split} #1 \end{split}\end{align*}}
\newcommand{\EQn}[1]{\begin{align}\begin{split} #1 \end{split}\end{align}}
\def\norm#1{\left\|#1\right\|}

\def\normb#1{\big\|#1\big\|}
\def\normbb#1{\Big\|#1\Big\|}

\def\brkb#1{\big(#1\big)}

\def\fbrk#1{\left\lbrace#1\right\rbrace}

\def\fbrkb#1{\big\lbrace#1\big\rbrace}

\def\wh#1{\widehat{#1}}
\def\wb#1{\overline{#1}}
\def\vec#1{\overrightarrow{#1}}

\def\pd{\partial}

\def\lsm{\lesssim}

\newcommand{\dd}{{\mathrm{d}}}

\def\al{\alpha}

\def\Om{\Omega}

\def\de{\delta}

\begin{document}

\setcounter{page}{1}

\title[Regularity of INLS]{Regularization for the inhomogeneous nonlinear Schr\"odinger equation in one dimension}

\author{Ruobing Bai}
\address{Ruobing Bai \newline School of Mathematics and Statistics\\
Henan University\\
Kaifeng, 475004, China \newline
Center for Applied Mathematics of Henan Province\\
Henan University\\
Zhengzhou, 450046, China
}
\email{baimaths@henu.edu.cn}
\thanks{}
\author{Luccas Campos}
\address{ Luccas Campos \newline  Department of Mathematics\\
Universidade Federal de Minas Gerais\\
Belo Horizonte,
MG31270-901
Brazil}
\email{luccas@ufmg.br}
\thanks{}
\author{Cui Ning}
\address{ Cui Ning \newline School of Financial Mathematics and Statistics\\
Guangdong University of Finance\\
Guangdong, 510521, China}
\email{cuining@gduf.edu.cn}
\thanks{}

%



\begin{abstract}\noindent
		In this work, we investigate the initial data problem for the following inhomogeneous nonlinear Schr\"odinger (INLS) equation
      $$  i \partial_t u + \pd_x^2 u = \pm |x|^{-b} |u|^{2k} u, \quad u(0, x)= u_0(x),$$
   where $b>0$ and $k\in \Z_{>0}$. The singular potential $|x|^{-b}$ impairs the Sobolev regularity that can be propagated by the flow. We establish the regularization mechanism of this model and prove local well-posedness in $H_x^{\alpha}$ for $\alpha=5/2-b-$. Moreover, at the endpoint space $H_x^{5/2-b}$, we show that there is no $C^{2k+1}$ flow map at the origin. The proof combines the one-dimensional local smoothing effect with a normal-form reduction for the high-low interaction generated by the singular part of the coefficient.  The argument applies to both the focusing and defocusing cases.
\end{abstract}

\maketitle

\section{Introduction}

Consider the Cauchy problem for the inhomogeneous nonlinear
Schr\"odinger equation
\begin{equation}\label{eq:NLS}\tag{INLS}
\left\{
\begin{aligned}
&i\partial_tu+\partial_x^2u=\lambda |x|^{-b}|u|^{2k}u,\\
&u(0,x)=u_0(x),
\end{aligned}
\right.
\end{equation}
where $u:\mathbb R_t\times\mathbb R_x\to\mathbb C$, $\lambda=\pm1$,
$b>0$, and $k\in\mathbb Z_{>0}$. The case $\lambda=1$ is referred to
as the defocusing case, while $\lambda=-1$ is referred to as the
focusing case. The spatially dependent coefficient $|x|^{-b}$ distinguishes
\eqref{eq:NLS} from the classical nonlinear Schr\"odinger equation and
introduces a singularity at the origin. This singularity is locally
integrable in dimension one when $0<b<1$, but its singular behavior imposes a
restriction on the Sobolev regularity that can be propagated by the
solution flow. The purpose of this paper is to determine this upper
regularity threshold in one space dimension.

The inhomogeneous nonlinear Schr\"odinger equation appears in models in
which the strength of the nonlinearity varies spatially. In particular,
such equations arise in the study of laser propagation through plasma
channels (see \cite{Gill-00,Liu-Tripathi-94}). From the mathematical
perspective, the equation has been extensively studied. In the
focusing case, Genoud and Stuart \cite{Genoud-Stuart-08} established local
well-posedness in the energy space $H^1(\mathbb R^N)$ for
$0<b<\min\{2,N\}$ and subcritical powers, together with global
well-posedness for arbitrary data in the mass-subcritical regime and for
small data in the intercritical regime. Farah \cite{Farah-16} later
obtained global well-posedness in $H^1(\mathbb R^N)$ in the focusing
intercritical case by means of a Gagliardo--Nirenberg inequality and a
quantitative smallness condition.

The local theory can be extended to fractional Sobolev spaces and fractional powers of the nonlinearity.
Let $N$ denote the spatial dimension, $\alpha\in \R$, $b,k \in \R_{>0}$ and define
\begin{equation}
k_\alpha^*
:=
\begin{cases}
\dfrac{2-b}{N-2\alpha},
&0\leq\alpha<\dfrac N2,\\[6pt]
\infty,
&\alpha\geq\dfrac N2.
\end{cases}
\end{equation}
For $0\leq\alpha<N/2$, the exponent $k=k_\alpha^*$ is the critical
power in $\dot H^\alpha(\mathbb R^N)$, while $k<k_\alpha^*$ corresponds
to the subcritical range. Guzm\'{a}n \cite{Guzman-17} used the contraction
mapping principle together with classical Strichartz estimates to prove
local and global well-posedness in $H^\alpha(\mathbb R^N)$ for subcritical parameters.
An and Kim \cite{A-K-21, A-K-23} and Aloui and Tayachi \cite{A-T-21} later treated the local theory for the critical and subcritical cases under different restrictions on
$b,\alpha,k$, and $N$.

More recently, Correia, Farah and the second author
\cite{C-C-F-25} obtained local well-posedness in
$H^\alpha(\mathbb R^N)$ with the minimal restrictions
$0<b<2$,
$0<b+\alpha<\min\left\{N,N/2+1\right\}$.
The question about well-posedness remained open in the range of parameters
\[
0<b<\min\left\{N,N/2+1\right\}\leq b+\alpha<N/2+2,
\]
since the integrals involved in the Picard iterations suggest convergence.
In dimension one, the range of parameters read
\[
0<b<1\leq b+\alpha<5/2.
\]
This is
the regime in which the singularity of the coefficient interacts directly
with the dispersive dynamics, and the convergence of the multilinear
expressions arising in the Picard iteration is no longer captured by the
standard estimates used at lower regularity.

The difficulty can be traced to the singular part of the coefficient.
Choose $\chi\in C_c^\infty(\mathbb R)$ such that $\chi=1$ in a
neighborhood of the origin and write
\begin{equation}\label{intro-potential-decomp}
|x|^{-b}=V+\widetilde V,
\qquad
V=\chi(x)|x|^{-b}.
\end{equation}
The remainder $\widetilde V$ is smooth with bounded derivatives of every
order, whereas the singular part satisfies
\begin{equation}\label{intro-V-regularity}
V\in L^q(\mathbb R)\cap H^\sigma(\mathbb R),
\qquad
1\leq q<\frac1b,
\qquad
\sigma<\frac12-b.
\end{equation}
Moreover,
\begin{equation}\label{intro-V-fourier}
|\widehat V(\xi)|\lesssim_b |\xi|^{b-1},
\qquad |\xi|\geq1.
\end{equation}
Thus $V$ is integrable throughout the range $0<b<1$, but its Fourier
transform has only algebraic decay. The distinction between the
integrability and the Sobolev regularity of $V$ is essential in the
high-regularity problem.

Indeed, if one applies the standard Strichartz argument directly to the
Duhamel term containing $V|u|^{2k}u$, then the desired Sobolev derivative
falls on the singular coefficient as well as on the solution. At the
regularity considered here, this procedure does not close. There is also
a spatial integrability obstruction: when $b\geq1/2$, the product
$V|u|^{2k}u$ need not belong to $L^2(\mathbb R)$ even for smooth $u$ that
does not vanish at the origin. Consequently, the singular term has to be
treated according to its frequency structure rather than as an ordinary
multiplication coefficient.

Writing
$G(u)=|u|^{2k}u$,
we decompose the singular interaction at a dyadic frequency $N$ as
\begin{equation}\label{intro-frequency-decomp}
P_N\bigl(VG(u)\bigr)
=
P_N\bigl(V(G(u)-G(u_{\ll N}))\bigr)
+
P_N\bigl(VG(u_{\ll N})\bigr).
\end{equation}
This separation distinguishes two different configurations.
In the first term, at least one factor of the solution carries frequency
comparable to the output frequency. In the second term, all solution
factors are much lower frequency and the high frequency is carried by the
singular coefficient itself.

The first configuration can be treated by transferring derivatives to a
high-frequency solution factor. The one-dimensional local smoothing
effect then provides the additional half derivative needed to close the
estimate. In particular, the dual local smoothing estimate
\begin{equation}\label{intro-smoothing}
\left\|
D^{1/2}
\int_0^t e^{i(t-s)\partial_x^2}F(s)\,ds
\right\|_{L_t^\infty L_x^2}
\lesssim
\|F\|_{L_x^1L_t^2}
\end{equation}
is compatible with the integrability $V\in L^1(\mathbb R)$ available for
every $0<b<1$. Thus the high-frequency solution factors provide a route
for the derivative to remain on the regular part of the dynamics.

The second configuration is more delicate. Here the singular coefficient
contains the high frequency, and there is no high-frequency solution
factor to absorb the Sobolev derivative. The favorable feature of this
interaction is instead temporal. Passing to the interaction representation
reveals an oscillatory phase
\begin{equation}\label{intro-phase}
\Phi(\vec\xi)
=
\left(\sum_{j=1}^{2k+2}\xi_j\right)^2
-
\sum_{j=2}^{2k+2}\iota_j\xi_j^2,
\end{equation}
where
\[
\iota_{2n}=1,
\qquad
\iota_{2n+1}=-1,
\qquad
n=1,\ldots,k,
\qquad
\iota_{2k+2}=1,
\]
and  $\xi_1$ denotes the frequency of $V$, while
$\xi_j$, $j=2,\ldots,2k+2$, denote the frequency of the solution
factors $u$ or $\overline u$. 

On the high-low region in which the output frequency and the frequency
carried by $V$ are comparable while all solution frequencies remain much
smaller, one has
\begin{equation}\label{intro-phase-size}
|\Phi(\vec\xi)|
\sim |\xi|^2
\sim |\xi_1|^2.
\end{equation}
Hence this interaction is temporally non-resonant. This non-resonance
allows us to integrate by parts in time and exploit the denominator
$\Phi^{-1}$.

The normal form method, introduced by Shatah \cite{Sha85CPAM}, is designed
precisely to exploit such non-resonant oscillations. Related regularization mechanisms for classical nonlinear Schr\"odinger equations with rough spatial potentials were studied by the first author, Lian, and Wu \cite{B-L-W-25, B-L-W-25-h} in one and higher dimensions using normal form methods. In the present
setting, the normal form reduction is adapted to the singular coefficient
rather than to the solution nonlinearity alone. The factor
$\Phi^{-1}$ contributes two powers of the large frequency in the
high-low interaction, effectively replacing the factor
$\langle\partial_x\rangle^\alpha V_H$ by an operator of order
$\langle\partial_x\rangle^{\alpha-2}V_H$. Schematically,
\begin{equation}\label{intro-normal-form}
\Omega\bigl(\langle\nabla\rangle^\alpha V_H,
u_L,\ldots,u_L\bigr)
\sim
\langle\nabla\rangle^{\alpha-2}V_H\,u_L\cdots u_L.
\end{equation}
The two derivatives gained from the oscillatory denominator are exactly
what is needed to compensate for the loss caused by the singular
coefficient.

At the dyadic level, \eqref{intro-V-fourier} gives
\[
\|P_NV\|_{L^2}\lesssim N^{b-\frac12},
\]
and therefore the transformed high-low term has size
\begin{equation}\label{intro-dyadic}
N^\alpha N^{-2}\|P_NV\|_{L^2}
\lesssim
N^{\alpha+b-\frac52}.
\end{equation}
This exhibits the exponent $5/2-b$ directly from the mechanism of
the proof. For
\[
\alpha=\frac52-b-\varepsilon,
\]
the right-hand side of \eqref{intro-dyadic} gains the summable factor
$N^{-\varepsilon}$. Consequently, the normal form transformation
recovers exactly the amount of regularity needed to close the
high-low interaction below the threshold.

The endpoint requires a different argument. The disappearance of the
factor $N^{-\varepsilon}$ in \eqref{intro-dyadic} could, in principle,
be a limitation of the method rather than a genuine obstruction of the
flow. To distinguish these possibilities, we examine the first nonlinear
Picard contribution for a suitable smooth initial datum whose Fourier
transform is supported in a fixed compact neighborhood of the origin.
Although all input factors are low frequency, the multiplication by the
singular coefficient generates a high-frequency output. On a suitable
family of disjoint high-frequency intervals, the corresponding Fourier
transform of
\[
A [u_0]=\int_0^t e^{-is\partial_x^2}(|x|^{-b} |e^{is\partial_x^2}u_0|^{2k}e^{is\partial_x^2}u_0)\dd s
\]
satisfies the lower bound
\begin{equation}\label{intro-endpoint-tail}
|\widehat{A[u_0]}(\xi)|
\gtrsim
|\xi|^{b-3}.
\end{equation}
Hence, at the endpoint
$$\alpha_*:=\frac52-b,$$
we have
\begin{equation}\label{intro-endpoint-density}
|\xi|^{2\alpha_*}
|\widehat{A[u_0]}(\xi)|^2
\gtrsim
|\xi|^{-1}.
\end{equation}
The resulting logarithmic divergence shows that the endpoint failure is
already encoded in the first nonlinear interaction. Thus the same
high-low mechanism has two complementary manifestations: below
$\alpha_*$, temporal non-resonance produces a two-derivative gain that
allows the flow to propagate the regularity, and at $\alpha_*$, the first
nonlinear iterate retains a non-summable high-frequency tail.

We therefore obtain the following sharp regularity statement.

\begin{thm}\label{thm-main}
Let $\lambda=\pm1$, $0<b<1$, and $k\in\mathbb Z_{>0}$. Then:

\begin{enumerate}
\item[(1)] for any $0<\varepsilon<(1-b)/2$ and
\[
\alpha=\frac52-b-\varepsilon,
\]
the Cauchy problem \eqref{eq:NLS} is locally well-posed in
$H^\alpha(\mathbb R)$;

\item[(2)] at the endpoint
\[
\alpha=\frac52-b,
\]
the data-to-solution map cannot be $C^{2k+1}$ at the origin as a map into
$C([0,T];H^\alpha(\mathbb R))$
on any nontrivial time interval $[0,T]$.
\end{enumerate}
\end{thm}

The theorem identifies $\alpha=\alpha_*$
as the sharp upper regularity threshold in the one-dimensional problem in
the sense that local well-posedness holds at every
$\alpha=\alpha_*-$, while the data-to-solution map fails to be
$C^{2k+1}$ at the endpoint. In particular, this extends the local
well-posedness range beyond the condition $b+\alpha<1$ covered in the
one-dimensional case of \cite{C-C-F-25}, and the endpoint obstruction
shows that the upper bound
$b+\alpha<5/2$
is not merely an artifact of the estimates.

The algebraic structure of the nonlinearity is essential in both parts of
the argument. Since
\[
|u|^{2k}u=u^{k+1}\overline{u}^{\,k},
\]
the nonlinear term admits a finite multilinear representation. This
structure is used in the normal form expansion, where the different
conjugated factors generate the phases in \eqref{intro-phase}, and also
in the $(2k+1)$-st derivative test for the data-to-solution map at the
endpoint. The argument applies to both $\lambda=1$ and $\lambda=-1$, as the sign of the nonlinear coefficient does not affect the non-resonant
structure of the high-low interaction or the endpoint divergence.

The rest of the paper is organized as follows. Section
\ref{sec:notations} introduces the notation, the decomposition
\eqref{intro-potential-decomp} of the singular coefficient, and the
high-low normal form transformation. Section \ref{sec-3} develops the
nonlinear estimates and proves Theorem \ref{thm-main}.

		\section{Preliminary}\label{sec:notations}
	
	\subsection{Notation} We use the following.

	$\bullet$ For any $a\in \mathbb{R}$, $a\pm:=a\pm\epsilon$ for arbitrarily small $\epsilon>0$.

$\bullet$ For any $z\in \mathbb{C}$, we define $\mbox{Re}\,z$ and $\mbox{Im}\,z$ as the real and imaginary parts of $z$, respectively.

$\bullet$ Denote $\langle \cdot \rangle=(1+|\cdot|^2)^{\frac12}$ and
$D^\alpha=(-\partial_x^2)^{\frac{\alpha}{2}}$.

$\bullet$ We write $X \lesssim Y$ or $Y \gtrsim X$ to indicate $X \leq CY$ for
	some constant $C>0$ independent of $X$ and $Y$. If $X \lesssim Y$ and $Y \lesssim X$, we write $X\sim Y$. If $X\leq \frac{2^{-5}}{2k+1} \, Y$, denote $X\ll Y$ or $Y\gg X$. Throughout the whole paper, the letter $C$ will denote a suitable positive constant that may vary from line to line.

$\bullet$	We also use the following norms to denote the mixed spaces $L_t^qL_x^r(I\times \R)$ and $L_x^r L_t^q( \R\times I)$, that is
	\begin{align*}
		\|u\|_{L_t^qL_x^r(I\times \R)}=\Big(\int_I \|u\|_{L_x^r( \R)}^q\dd t\Big)^{\frac{1}{q}}, \quad \|u\|_{L_x^rL_t^q(\R\times I)}
&=\Big(\int_\R \|u\|_{L_t^q(I)}^r
\dd x\Big)^{1/r}.
	\end{align*}
For simplicity, we often write $L_t^{q}L_x^r:=L_t^{q}L_x^r(I\times \R)$, $L_x^rL_t^{q}:=L_x^rL_t^{q}( \R\times I)$ and some similar simplified norm notations for short.

$\bullet$ $\widehat{f}$ or $\mathscr{F}f$ denotes the Fourier transform of $f$:
$$
\mathscr{F}f (\xi)= \widehat f(\xi)
	= \int_{\R} e^{- i x\cdot \xi} f(x) \dd x.
$$
We also define
$$
\mathscr{F}^{-1}g (x)= (2\pi)^{-1}\int_{\R} e^{ i x\cdot \xi} g(\xi) \dd \xi.
$$

$\bullet$	We fix a nonnegative cut-off function $\phi(r) \in C_c^{\infty}{([0, +\infty))}$ such that $\phi(r)=1$ if $r\leq 1$ and $\phi(r)=0$ if $r>2$.
For dyadic number $N\in 2^{\mathbb{N}}$, when $N\geq 1$, let $\phi_{\leq N}(r)=\phi(N^{-1}r)$. Then, we define
$\phi_{1}(r):=\phi(r)$, and $\phi_{N}(r)=\phi_{\leq N}(r)-\phi_{\leq \frac N2}(r)$ for any $N\geq 2$. We define the inhomogeneous Littlewood-Paley dyadic operator
$$
f_1=P_{ 1}f:=\mathscr{F}^{-1}(\phi_{ 1}(|\xi|)\widehat{f}(\xi)),
$$
and for any $N\geq 2$,
$$
f_{ N}=P_{ N}f:=\mathscr{F}^{-1}(\phi_{ N}(|\xi|)\widehat{f}(\xi)).
$$
By definition, we have $f=\sum_{N\in 2^{\N}}f_N$. Moreover, we also define the following:
	\begin{align*}
	&f_{\leq N}=P_{\leq N}f:=\mathscr{F}^{-1}(\phi_{\leq N}(|\xi|)\widehat{f}(\xi)),\\
	&f_{ \ll N}=P_{ \ll N}f:=\mathscr{F}^{-1}(\phi_{ \leq N}(2^5(2k+1)|\xi|)\widehat{f}(\xi)),\\
&f_{\lsm N}=P_{\lsm N}f:=\mathscr{F}^{-1}(\phi_{ \leq N}(2^{-5}{(2k+1)^{-1}}|\xi|)\widehat{f}(\xi)).
	\end{align*}
	We also define $f_{\geq N}=P_{\geq N}f:=f-f_{\leq N}$, $f_{\gg N}=P_{\gg N}f:=f-P_{\lsm N}f$, and $f_{\gtrsim N}=P_{\gtrsim N}f:=f-P_{ \ll N}f$.


	$\bullet$ For any $1\leq p <\infty$, we define $l_{N}^p=l_{N\in 2^{\mathbb{N}}}^p$ by its norm
	$$
	\|c_N\|_{l_{N\in 2^{\mathbb{N}}}^p}^p:=\sum_{N\in 2^{\mathbb{N}}}|c_N|^p.
	$$
	For $p =\infty$, we define $l_{N}^{\infty}=l_{N\in 2^{\mathbb{N}}}^{\infty}$ by its norm
	$$
	\|c_N\|_{l_{N\in 2^{\mathbb{N}}}^{\infty}}:=\sup_{N\in 2^{\mathbb{N}}}|c_N|.
	$$

	$\bullet$ We define the Triebel-Lizorkin spaces $F_{p}^{\alpha, q}$ with the corresponding norm as follows,
	\EQ{
		\|u\|_{F_{p}^{\alpha, q}}=\|u\|_{L_x^p}+\|N^{\alpha}P_Nu\|_{ L_x^p l_{N\in 2^{\mathbb{N}}}^q}.
	}

$\bullet$
If the pair $(q, r)\in \R^2$ satisfies
		\begin{align*}
			2\leq q, r\leq \infty  ,\quad \frac{2}{q}+\frac{1}{r}=\frac{1}{2},
		\end{align*}
we say that the pair $(q, r)$ is \textit{$L^2$-admissible}.

	\subsection{Useful lemmas}

	\quad In this section, we state some preliminary estimates that will be used in our later sections.
	Firstly, we introduce the following Bernstein estimates.
\begin{lem}[Bernstein estimates\label{lem:Bernstein}]
		For any $1\leq p \leq q \leq \infty$, $s\geq 0$, and $f\in L_x^p(\R^d)$,
		\begin{align*}
			\||\nabla|^{\pm s}P_N f\|_{L_x^p(\R^d)}&\sim N^{\pm s}\|P_N f\|_{L_x^p(\R^d)},\quad N>1\\
			\|P_{\leq N} f\|_{L_x^q(\R^d)}&\lesssim  N^{\frac dp-\frac dq}\|P_{\leq N} f\|_{L_x^p(\R^d)},\\
			\|P_N f\|_{L_x^q(\R^d)}&\lesssim  N^{\frac dp-\frac dq}\|P_N f\|_{L_x^p(\R^d)}.
		\end{align*}
	\end{lem}

	\begin{lem}[Schur's test\label{lem:Schur}]
		For any $a>0$, let sequences $\{a_N\}, \{b_N\}\in l_{N\in 2^{\N}}^2$, then we have
		\EQ{
			\sum_{ N\geq N_1}\Big({\frac{N_1}{N}}\Big)^a a_N b_{N_1}\lesssim\|a_N\|_{l_N^2}\|b_N\|_{l_N^2}.
		}
	\end{lem}

	Next, we recall the well-known Strichartz estimates.
	\begin{lem}
		[Strichartz estimates, see \cite{Keel-Tao-1998, Cazenave-03, Strichartz-1977, Ginibre-Velo-1985} \label{lem:strichartz}] Let $I\subset \R$ be a time interval. Let $(q_j, r_j), j=1,2,$ be $L^2$-admissible,
		then the following statements hold:
		\begin{align}\label{1.2222}
			\|e^{it\pd_x^2}f\|_{L_t^{q_j}L_x^{r_j}(I\times{\R})}\lesssim\|f\|_{L^2(\R)};
		\end{align}
		and
		\begin{align}\label{1.222234}
			\Big\|\int_0^t e^{i(t-s)\pd_x^2}F(s)\dd s\Big\|_{L_t^{q_1}L_x^{r_1}(I\times{\R})}\lesssim \|F\|_{L_t^{q_2'}L_x^{r_2'}(I\times{\R})},
		\end{align}
		where $\frac{1}{q_2}+\frac{1}{q_2'}=\frac{1}{r_2}+\frac{1}{r_2'}=1$.
	\end{lem}

The next lemma is the smoothing effects.
	\begin{lem}
		[Smoothing effects, see \cite{KenigPonceVega-CPAM-93, LinaresPonce-09}] Let $I\subset \R$ be an interval, including $I=\R$. Then\\
		1)
		\EQn{\label{Smooth1}
			\|D^{\frac 12}e^{it\pd_x^2}f \|_{L_x^\infty L_t^2(\R\times I )} \lesssim  \| f  \|_{L_x^2(\R)},
		}
		for all $f\in L_x^2(\R)$; and\\
		2)
		\EQn{\label{Smooth2}
			\Big\|  D^{\frac 12} \int_0^t e^{i(t-s)\pd_x^2} F(x,s)\dd s\Big\|_{L_t^\infty L_x^2(I\times \R)}
			\lesssim           \|F\|_{L_x^1 L_t^2(\R\times I)}.
		}
	
	\end{lem}

%
%
%

			We also need the following Littlewood-Paley theory, see Remark 2.2.2 in \cite{Gra-14}.
	\begin{lem}[Littlewood-Paley theory\label{lem:littlewood-Paley}]
		Let $1<p<\infty$, for any $\alpha>0$, we have
		\EQ{
			\|f\|_{F_{p}^{\alpha, 2}}\sim \|\langle\nabla\rangle^{\alpha}f\|_{L_x^p}.
		}
	\end{lem}

	Next, we recall the Coifman-Meyer multiplier theory.
	
	\begin{lem}[Multilinear Coifman-Meyer multiplier estimates, see \cite{Co-Me-91}\label{lem:Coifman-Meyer}]
		Let the function $m$ on ${(\R)^k}$ be bounded and let $T_m$ be the corresponding $k$-linear multiplier operator on $\R$
		\begin{align*}
			T_m(f_1,\cdots, f_k)(x)=\int_{\R^k}m(\xi_1,\cdots,\xi_k)\hat{f_1}(\xi_1)\cdots\hat{f_k}(\xi_k)e^{ix\cdot(\xi_1+\cdots+\xi_k)}d\xi_1\cdots d\xi_k.
		\end{align*}
		If $L$ is sufficiently large and $m$ satisfies
		\begin{align}\label{condition-m}
			\Big|\partial_{\xi_1}^{\al_1}\cdots\partial_{\xi_k}^{\al_k}m(\xi_1,\cdots,\xi_k)\Big|\lesssim_{\al_1,\cdots, \al_k}(|\xi_1|+\cdots+|\xi_k|)^{-(|\al_1|+\cdots+|\al_k|)},
		\end{align}
		for multi-indices $\al_1,\cdots, \al_k$ satisfying $|\al_1|+\cdots+|\al_k|\leq L$. Then, for $1< p<\infty$, $1<p_1, \cdots, p_k\leq \infty$ and $\frac 1p=\frac{1}{p_1}+\cdots+\frac{1}{p_k}$, we have
		\begin{align*}
			\|T_m(f_1,\cdots, f_k)\|_{L_x^p}\leq C\|f_1\|_{L_x^{p_1}}\cdots\|f_k\|_{L_x^{p_k}}.
		\end{align*}
	\end{lem}
The Coifman-Meyer Multiplier Theorem is reduced to the Mihlin-H\"{o}rmander Multiplier
Theorem when $k=1$ and $1<p<\infty$.

We also need the following vector valued multiplier estimates, see Lemma 2.9 in \cite{B-L-W-25}.
\begin{lem}\label{vector-estimate}
Let the function $m$ on ${(\R^d)^k}$ be defined as in Lemma \ref{lem:Coifman-Meyer}. Let $T_m$ be the $k$-linear multiplier operator on $\R\times \R^d$, where $\R^d$ is the $d$-dimensional Euclidean space,
		\begin{align*}
			T_m(f_1,\cdots, f_k)(t, x)=\int_{(\R^d)^k}m(\xi_1,\cdots,\xi_k)\hat{f_1}(t, \xi_1)\cdots\hat{f_k}(t, \xi_k)e^{ix\cdot(\xi_1+\cdots+\xi_k)}d\xi_1\cdots d\xi_k.
		\end{align*}	
		Then, for $1< p, q<\infty$, $1<p_1, \cdots, p_k, q_1, \cdots, q_k < \infty$ and $\frac 1p=\frac{1}{p_1}+\cdots+\frac{1}{p_k}$, $\frac 1q=\frac{1}{q_1}+\cdots+\frac{1}{q_k}$, we have
		\begin{align*}
			\|T_m(f_1,\cdots, f_k)\|_{L_x^pL_t^q}\leq C\|f_1\|_{L_x^{p_1}L_t^{q_1}}\cdots\|f_k\|_{L_x^{p_k}L_t^{q_k}}.
		\end{align*}
In particular, if $f_1, \cdots, f_j$ $(j< k)$ are independent of time variable $t$, then for $1<p, p_1, \cdots, p_k, q, \\
q_{j+1}, \cdots, q_k < \infty$ and $\frac 1p=\frac{1}{p_1}+\cdots+\frac{1}{p_k}$, $\frac 1q=\frac{1}{q_{j+1}}+\cdots+\frac{1}{q_k}$,
\begin{align}\label{Coifman-meyer-vector}
			\|T_m(f_1,\cdots, f_k)\|_{L_x^pL_t^q}\leq C\|f_1\|_{L_x^{p_1}}\cdots\|f_j\|_{L_x^{p_j}}\|f_{j+1}\|_{L_x^{p_{j+1}}L_t^{q_{j+1}}}\cdots\|f_k\|_{L_x^{p_k}L_t^{q_k}}.
		\end{align}
	\end{lem}
The next result is the endpoint Kato-Ponce inequality, which was established in \cite{Bou-Li-12, Li-19} and originally proved in \cite{Ka-Po-88}.
\begin{lem}[Kato-Ponce inequality]
For $s>0$, $1<p<\infty$, and
$1<p_1,p_2,p_3,p_4\leq\infty$ satisfying
\[
\frac{1}{p}=\frac{1}{p_1}+\frac{1}{p_2},
\qquad
\frac{1}{p}=\frac{1}{p_3}+\frac{1}{p_4},
\]
the following Kato--Ponce inequality holds:
\begin{equation}
\|D^s(fg)\|_{L^p(\mathbb{R})}
\lesssim
\|D^s f\|_{L^{p_1}(\mathbb{R})}
\|g\|_{L^{p_2}(\mathbb{R})}
+
\|D^s g\|_{L^{p_3}(\mathbb{R})}
\|f\|_{L^{p_4}(\mathbb{R})}.
\end{equation}
\end{lem}

\subsection{Normal form transform}
Next, we apply the normal form transform to the equation \eqref{eq:NLS} with $p=2k$, and $k\in \N$.
Then, the integral equation for \eqref{eq:NLS} is
\EQn{\label{eq:integral}
u(t,x) = e^{it\pd_x^2} u_0(x) - i\lambda\int_{0}^t e^{i(t-s)\pd_x^2}|x|^{-b}|u|^{2k}u(s)\dd s.
}
We take a cut-off function $\chi(x) \in C_c^{\infty}(\R)$ such that
	$$
	\chi(x) =\left\{ \aligned
	&1, \quad |x| \leq 1,
	\\
	&0, \quad |x|\geq 2.
	\endaligned
	\right.
	$$
 We define
\begin{align}\label{VV}
V(x)\triangleq \chi(x)|x|^{-b} \mbox{ and }  \widetilde{V}(x)\triangleq (1-\chi(x))|x|^{-b}.
\end{align}
 Then the integral equation \eqref{eq:integral} reads as
\EQn{\label{eq:integral-2}
u(t,x) = e^{it\pd_x^2} u_0(x) -& i\lambda\int_{0}^t e^{i(t-s)\pd_x^2}V(x)|u|^{2k}u(s)\dd s\\
-&i\lambda\int_{0}^t e^{i(t-s)\pd_x^2}\widetilde{V}(x)|u|^{2k}u(s)\dd s.
}
The following elementary result is the key point that allows us to close the nonlinear estimates in the full range $0<b<1$.
\begin{lem}\label{lem:potential-regularity}
Let $0<b<1$, and let $V$ and $\widetilde V$ be given by \eqref{VV}. Then, for every $1\leq q<1/b$ and every $\sigma<\frac12-b$,
\begin{equation}\label{V-regularity}
V\in  L^{q}(\R)\cap H^\sigma(\R),
\end{equation}
and for every integer $m\geq0$, $\widetilde V\in W^{m,\infty}(\R)$. More generally, for every real $m\geq0$, $\langle\partial_x\rangle^m\widetilde V\in L^\infty(\R)$.
\end{lem}

\begin{proof}
Since $V$ is supported in $[-2,2]$, $0<b<1$, and $1\leq q<1/b$,
$$
\|V\|_{L^q}^q
\leq 2\int_0^2x^{-bq}\dd x<\infty.
$$
Next, we estimate $\|V\|_{H^\sigma(\R)}$. We consider two cases separately: $|\xi|\geq 1$ and $|\xi|\leq 1$.
Fix $|\xi|\geq 1$, and set $ r=|\xi|^{-1}$.  Decompose $\widehat V(\xi)$ as follows
$$
\widehat V(\xi)
=\int_{|x|\leq r}e^{-ix\xi}V(x)\dd x
 +\int_{r\leq|x|\leq2}e^{-ix\xi}V(x)\dd x
=:J_0(\xi)+J_1(\xi).
$$
For the first part $J_0(\xi)$,
\begin{equation}\label{V-small-x}
|J_0(\xi)|
\lsm \int_0^r x^{-b}\,\dd x
=\frac{1}{1-b}|\xi|^{b-1}.
\end{equation}
For the second part $J_1(\xi)$, treat the positive and negative half-lines
separately.  On $[r,2]$, let $a(x)=\chi(x)x^{-b}$.  Integration-by-parts gives
\begin{align*}
\left|\int_r^2e^{-ix\xi}a(x)\dd x\right|
&\leq |\xi|^{-1}\big(|a(r)|+|a(2)|\big)
   +|\xi|^{-1}\int_r^2|a'(x)|\dd x\\
&\lesssim |\xi|^{-1}r^{-b}
 +|\xi|^{-1}\int_r^2\big(x^{-b-1}+x^{-b}\big)\dd x\\
&\lesssim |\xi|^{b-1}.
\end{align*}
The negative half-line $[-2, -r]$ can be treated similarly after the change of variables
$x\mapsto-x $.  Together with \eqref{V-small-x}, this proves
$$
|\widehat V(\xi)|\lesssim_b|\xi|^{b-1},
\qquad |\xi|\geq1.
$$
For $|\xi|\leq 1$, one has
$|\widehat V(\xi)|\leq\|V\|_{L^1}$.  Hence, we conclude that
\begin{align*}
\|V\|_{H^\sigma}^2
&=\int_{|\xi|\leq 1}\langle\xi\rangle^{2\sigma}
 |\widehat V(\xi)|^2\dd\xi
 +\int_{|\xi|\geq 1}\langle\xi\rangle^{2\sigma}
 |\widehat V(\xi)|^2\dd\xi\\
&\lesssim \|V\|_{L^1}^2
 +\int_1^\infty \xi^{2\sigma+2b-2}\dd\xi.
\end{align*}
The last integral is finite if and only if
$2\sigma+2b-2<-1$, namely
$\sigma<\frac12-b$. This proves \eqref{V-regularity}.

Next, we prove the asserted $L^\infty$ bounds. By the definition of $\widetilde V$, we directly obtain that $\widetilde V\in W^{j,\infty}(\R)$ for every integer $j\geq0$. To extend this conclusion to arbitrary real numbers $m\geq0$, choose an integer $M>m$. Since $\langle\partial_x\rangle^m=\langle\partial_x\rangle^{-(M-m)}\langle\partial_x\rangle^M$, it suffices to prove that $\langle\partial_x\rangle^{-\delta}$ is bounded on $L^\infty$ for every $\delta>0$ and use the integer-order bounds for $\langle\partial_x\rangle^M\widetilde V$, that is $\|\langle\pd_x \rangle^{-\de} \widetilde V\|_{L_x^{\infty}}\lsm \| \widetilde V\|_{L_x^{\infty}}$, for any $\de>0$. Indeed, we firstly have that for any $\de>0$,
$$
\langle\pd_x\rangle^{-\de} \widetilde V=\mathscr{F}^{-1}\big((1+|\xi|^2)^{-\frac {\de} 2}\big)\ast \widetilde V.
$$
Notice that
$$
(1+|\xi|^2)^{-\frac {\de} 2}=\frac 1{\Gamma(\frac {\de} 2)}\int_0^{\infty}e^{-t(1+|\xi|^2)}t^{\frac {\de}{2}-1}dt,
$$
where the Gamma function $\Gamma(s)=\int_0^{\infty}e^{-t}t^{s-1}dt$. Furthermore, $\mathscr{F}^{-1}( e^{-t|\xi|^2})=(2\sqrt \pi)^{-1}t^{-\frac 12}e^{-\frac {|x|^2} {4t}}$.
Hence, we conclude that
$$
\mathscr{F}^{-1}\big((1+|\xi|^2)^{-\frac {\de} 2}\big)=\frac 1{2\sqrt \pi\,\Gamma(\frac {\de} 2)}\int_0^{\infty}e^{-t}e^{-\frac {|x|^2} {4t}}t^{\frac {\de}{2}-\frac 32}dt.
$$
Further, we obtain
$$
\|\mathscr{F}^{-1}\big((1+|\xi|^2)^{-\frac {\de} 2}\big)\|_{L_x^1}\lsm \frac 1{\Gamma(\frac {\de} 2)}\int_0^{\infty}e^{-t}t^{\frac 12 }t^{\frac {\de}{2}-\frac 32}dt=1.
$$
Therefore, by the Young's inequality for convolutions,
$$
\|\langle\pd_x\rangle^{-\de} \widetilde V\|_{L_x^{\infty}}\lsm \|\mathscr{F}^{-1}\big((1+|\xi|^2)^{-\frac {\de} 2}\big)\|_{L_x^{1}}\|\widetilde V\|_{L_x^{\infty}}\lsm \|\widetilde V\|_{L_x^{\infty}}.
$$
This proves \eqref{V-regularity}, and thus completes the proof of this lemma.
\end{proof}

Next, we define the normal form transform with respect to some suitable frequency part of the nonlinear term $V(x)|u|^{2k}u$. More precisely, the target regularity is $H_x^{2+\frac 12-b-}(\R)$ while $V(x) \in H_x^{\frac 12-b-}$. Hence, when the factor $V$ has high frequency while the other functions $u$ remain low-frequency, this scenario becomes the most difficult to address. In this case, observe that the phase function does not vanish, allowing us to define the normal form transformation for this term.
\begin{definition}[Normal form transform]
\label{defn:normal-form}
Let $k\in\Z$ with $k\ge1$ and $\xi,\xi_j\in\R$ for $j=1,2,3,\cdots,2k+2$. Denote that $\vec \xi=(\xi_1,\xi_2,\xi_3,\cdots,\xi_{2k+2})\in\R^{2k+2}$. Let $\iota_j\in\fbrk{+1,-1}$ for $j=2,3,\cdots,2k+1, 2k+2$ such that
\EQ{
\iota_{2n} := +1\text{, } \iota_{2n+1} :=-1 \text{, for }n=1,2,\cdots,k,
}
and
\EQ{
\iota_{2k+2}=+1.
}
Define the hypersurface
\EQ{
\Ga(\vec \xi) := \fbrkb{\vec \xi\in\R^{2k+2}:\xi=\sum_{j=1}^{2k+2}\xi_j}.
}
The phase function is defined by
\EQ{
\Phi(\vec \xi):=(\sum_{j=1}^{2k+2}\xi_j)^2- \sum_{j=2}^{2k+2}\iota_j \xi_j^2.
}
Denote the multiplier
\EQ{
m_N(\vec \xi):=\frac {\langle\xi\rangle^{\alpha}\langle\xi_1\rangle^{2-\alpha}}{\Phi( \vec \xi)} \prod_{j=2}^{2k+2} \phi_{\ll N}(|\xi_j|)\phi_{N}(|\xi|).
}
For any function $f$, we define
\EQ{
f^{(+1)} := f\text{, and }f^{(-1)} := \wb f.
}
Using the above notations, we give the following definitions:
\begin{enumerate}
\item (Normal form transform)
Let $N_0\in 2^\N$ be sufficiently large. We define the normal form transform for functions $f_1,f_2,\cdots,f_{2k+2}$ by
\EQ{
\Om(f_1,f_2,\cdots,f_{2k+2})(x) := \sum_{N:N\geq N_0}\int_{\Ga(\vec\xi)} e^{ix\brkb{\sum_{j=1}^{2k+2}\xi_j}} m_N(\vec \xi) \wh{f_1}(\xi_1)\prod_{j=2}^{2k+2} \wh{f_j^{(\iota_j)}}(\xi_j) \dd \vec\xi.
}
\item (Simplified notation for normal form)
If $f_2=\cdots=f_{2k+2}=f$, we denote that
\EQ{
	\Om(f_1, f):=\Om(f_1,f_2,\cdots,f_{2k+2})=\Om(f_1,f,\cdots,f).
}
If $f_j=f$ for all $2\le j\le 2k+2$ except for $j=l$, and $f_l=g$, then we denote that
\EQ{
	\Om_l(f_1, f, g):=\Om(f_1,f_2,\cdots,f_{2k+2})=\Om(f_1,f,\ldots,f,\underbrace{g}_{\ell\text{-th slot}},f,\ldots,f)
}

\end{enumerate}

\end{definition}

The definition yields the following identity.
\begin{lem} \label{lem:decomposition}
Let $u(t,x )$ solve \eqref{eq:integral-2}, $V$ and $\widetilde{V}$ be defined in \eqref{VV}, the multi-linear operator $\Om$ and $\Om_l$ be defined in Definition \ref{defn:normal-form}. Then for any $\alpha \in \R$, with
$$
F(V,\widetilde V,u):=(V+\widetilde V)|u|^{2k}u,
$$
the following identity holds:
\begin{align*}
\langle\partial_x\rangle^{\alpha} u(t)
=&\langle\partial_x\rangle^{\alpha}e^{i t\pd_x^2}u_0-i\lambda\int_{0}^t \langle\partial_x\rangle^{\alpha} e^{i(t-s)\pd_x^2}\widetilde{V}(x)|u|^{2k}u(s)\dd s\\
&-i\lambda\sum_{N\in 2^{\N}}\int_0^t \langle\partial_x\rangle^{\alpha}e^{i(t-s)\pd_x^2}P_N\big(V(|u|^{2k}u-|u_{\ll N}|^{2k}u_{\ll N})\big)\dd s\notag \\
    &-i\lambda\sum_{N: N< N_0}\int_0^t \langle\partial_x\rangle^{\alpha}e^{i(t-s)\pd_x^2}P_N\big(V|u_{\ll N}|^{2k}u_{\ll N}\big)\dd s\\
    &-\lambda \Om(\langle\partial_x\rangle^{\alpha-2}V, u)+\lambda e^{it\pd_x^2}\Om(\langle\partial_x\rangle^{\alpha-2}V, u_0)\notag\\
&+i\lambda^2 \sum_{l=2}^{2k+2}(-1)^{l-1}\int_0^t e^{i(t-s)\pd_x^2}\Om_l(\langle\partial_x\rangle^{\alpha-2}V, u, F(V, \widetilde{V}, u))(s)\dd s.
    \end{align*}
\end{lem}
\begin{proof}
By \eqref{eq:integral-2} and the high-low frequency decomposition, we have
\begin{align}\label{It}
\langle\partial_x\rangle^{\alpha}  u(t)
  =&\langle\partial_x\rangle^{\alpha}e^{i t\pd_x^2}u_0 -i\lambda\int_{0}^t \langle\partial_x\rangle^{\alpha} e^{i(t-s)\pd_x^2}\widetilde{V}(x)|u|^{2k}u(s)\dd s\notag\\
  &-i\lambda\sum_{N\in 2^{\N}}\int_0^t \langle\partial_x\rangle^{\alpha}e^{i(t-s)\pd_x^2}P_N\big(V(|u|^{2k}u-|u_{\ll N}|^{2k}u_{\ll N})\big)\dd s\notag \\
    &-i\lambda\sum_{N: N< N_0}\int_0^t \langle\partial_x\rangle^{\alpha}e^{i(t-s)\pd_x^2}P_N\big(V|u_{\ll N}|^{2k}u_{\ll N}\big)\dd s
    + I(t),
\end{align}
where
\begin{align*}
  I(t)
  :=-i\lambda\sum_{N: N\geq N_0}\int_0^t \langle\partial_x\rangle^{\alpha} e^{i(t-s)\pd_x^2}P_N\big(V|u_{\ll N}|^{2k}u_{\ll N}\big)\dd s .
\end{align*}
Let $u=e^{it\partial_x^2}v$, using the Fourier transform, we have
\begin{align*}
  \widehat{I}(\xi)=&-i\lambda\sum_{N: N\geq N_0}\int_0^t\int_{\Ga(\vec\xi)}e^{-it|\xi|^2}  e^{is\Phi(\vec \xi)}\langle\xi\rangle^{\alpha}\phi_{N}(|\xi|)\widehat{V}(\xi_1)\\
  &\qquad\qquad \cdot\prod_{j=2}^{2k+2} \phi_{\ll N}(|\xi_j|) \wh{v^{(\iota_j)}}(\xi_j)\dd \xi_1\cdots \dd \xi_{2k+1}\dd s,
  \end{align*}
  where
\EQ{
\Phi( \vec \xi):=|\xi|^2- \sum_{j=2}^{2k+2}\iota_j \xi_j^2\sim |\xi|^2\sim |\xi_1|^2.
}
Since this term is non-resonant, we can use the integration-by-parts to treat it. Here, we
note that
\EQ{
			\partial_{s}\widehat{v}(\xi)=-i\lambda e^{is|\xi|^2}\widehat{F(V, \widetilde{V}, u)}(\xi).
		}

{

Hence, we have
\begin{align*}
\widehat{I}(\xi)=&-i\lambda\sum_{N: N\geq N_0} e^{-it|\xi|^2}
\int_{\Ga(\vec\xi)}
\frac {e^{it \Phi( \vec \xi) }}{i\Phi( \vec \xi)}
\langle\xi\rangle^{\alpha}
\phi_{N}(|\xi|)
\prod_{j=2}^{2k+2} \phi_{\ll N}(|\xi_j|)
\widehat{V}(\xi_1)
\wh{v^{(\iota_j)}}(\xi_j)
\dd \xi_1\cdots \dd \xi_{2k+1}\\
&+i\lambda\sum_{N: N\geq N_0} e^{-it|\xi|^2}
\int_{\Ga(\vec\xi)}
\frac {1}{i\Phi( \vec \xi)}
\langle\xi\rangle^{\alpha}
\phi_{N}(|\xi|)
\prod_{j=2}^{2k+2} \phi_{\ll N}(|\xi_j|)
\widehat{V}(\xi_1)
\wh{u_0^{(\iota_j)}}(\xi_j)
\dd \xi_1\cdots \dd \xi_{2k+1}\\
&-\lambda^2 \sum_{N: N\geq N_0}
\sum_{l=2}^{2k+2}(-1)^{l-1}
\int_0^t\int_{\Ga(\vec\xi)}
e^{-it|\xi|^2}
\frac {e^{is \Phi( \vec \xi) }}{i\Phi( \vec \xi)}
\langle\xi\rangle^{\alpha}
\phi_{N}(|\xi|)
\prod_{j=2}^{2k+2} \phi_{\ll N}(|\xi_j|)
\widehat{V}(\xi_1)
e^{i{\color{red}\iota_l}s|\xi_l|^2}\\
&\cdot\wh{v^{(\iota_2)}}(\xi_2)\cdots
\wh{v^{(\iota_{l-1})}}(\xi_{l-1})
\wh{F^{(\iota_l)}(V, \widetilde{V}, u)}(\xi_l)
\wh{v^{(\iota_{l+1})}}(\xi_{l+1})\cdots
\wh{v^{(\iota_{2k+2})}}(\xi_{2k+2})
\dd \xi_1\cdots \dd \xi_{2k+1} \dd s.
\end{align*}

Using the notation of the multiplier $m_N(\vec \xi)$, we rewrite the above identity as follows,
\begin{align*}
\widehat{I}(\xi)=&-\lambda\sum_{N: N\geq N_0}
\int_{\Ga(\vec\xi)}
m_N(\vec \xi)
\langle\xi_1\rangle^{\alpha-2}
\widehat{V}(\xi_1)
\prod_{j=2}^{2k+2}
\wh{u^{(\iota_j)}}(\xi_j)
\dd \xi_1\cdots \dd \xi_{2k+1}\\
&+\lambda\sum_{N: N\geq N_0}
e^{-it|\xi|^2}
\int_{\Ga(\vec\xi)}
m_N(\vec \xi)
\langle\xi_1\rangle^{\alpha-2}
\widehat{V}(\xi_1)
\prod_{j=2}^{2k+2}
\wh{u_0^{(\iota_j)}}(\xi_j)
\dd \xi_1\cdots \dd \xi_{2k+1}\\
&+i\lambda^2 \sum_{N: N\geq N_0}
\sum_{l=2}^{2k+2}(-1)^{l-1}
\int_0^t\int_{\Ga(\vec\xi)}
e^{-i(t-s)|\xi|^2}
m_N(\vec \xi)
\langle\xi_1\rangle^{\alpha-2}
\widehat{V}(\xi_1)\\
&\cdot\wh{u^{(\iota_2)}}(\xi_2)\cdots
\wh{u^{(\iota_{l-1})}}(\xi_{l-1})
\wh{F^{(\iota_l)}(V, \widetilde{V}, u)}(\xi_l)
\wh{u^{(\iota_{l+1})}}(\xi_{l+1})\cdots
\wh{u^{(\iota_{2k+2})}}(\xi_{2k+2})
\dd \xi_1\cdots \dd \xi_{2k+1} \dd s.
\end{align*}
}
Further, by the definition of multi-linear operator $\Om$ and $\Om_l$, we have
\begin{align}\label{I1}
I={}&-\lambda \Om(\langle\partial_x\rangle^{\alpha-2}V, u)
+\lambda e^{it\pd_x^2}
\Om(\langle\partial_x\rangle^{\alpha-2}V, u_0)\notag\\
&+i\lambda^2
\sum_{l=2}^{2k+2}(-1)^{l-1}
\int_0^t
e^{i(t-s)\pd_x^2}
\Om_l(\langle\partial_x\rangle^{\alpha-2}V,
u,F(V, \widetilde{V}, u))(s)\dd s.
\end{align}
Combining \eqref{It} and \eqref{I1}, this completes the proof of this lemma.
\end{proof}

\section{Proof of Theorem \ref{thm-main}} \label{sec-3}

Throughout this section,
let $\epsilon$ be a fixed sufficiently small constant satisfying $0<2\epsilon<1-b$. Let
$$\alpha=\frac 52-b-\epsilon.$$ Note that $\alpha> 3/2$.
 Next, we provide the necessary nonlinear estimates.
\subsection{Easier nonlinear terms}
The crucial point is that the following estimate uses the
$L^1$-integrability of $V$. For convenience, let
$$
G(u):=|u|^{2k}u.
$$
\begin{lem}\label{lem:easy}
		Let $I=[0, T) \subset\R^+$ with $T<1$. Then for any $N_0\in 2^{\N}$,
		\EQn{\label{easy-Vuu}
		\normbb{\int_{0}^te^{i(t-s)\pd_x^2}\langle\partial_x\rangle^{\alpha} \big(\widetilde{V}G(u)\big)(s)\dd s}_{L_t^{\infty}L_x^2(I\times \R)} \lsm T\|\langle\partial_x\rangle^{\alpha} \widetilde{V}\|_{L_x^{\infty}}\|u\|_{L_{t}^{\infty}H_x^{\alpha}}^{2k+1},
		}
\EQn{\label{easy-Vuu2}
		\normbb{\sum_{N\in 2^{\N}}\int_0^t \langle\partial_x\rangle^{\alpha}e^{i(t-s)\pd_x^2}P_N\big(V(G(u)-G(u_{\ll N}))\big)(s)\dd s}_{L_t^{\infty}L_x^2(I\times \R)} \lsm T^{\frac 12}\|V\|_{L_x^1}\|u\|_{L_t^{\infty}H_x^{\alpha}}^{2k+1},
		}
and
\EQn{\label{easy-Vuu3}
		\normbb{\sum_{N: N< N_0}\int_0^t \langle\partial_x\rangle^{\alpha}e^{i(t-s)\pd_x^2}P_N\big(V G(u_{\ll N} )\big)(s)\dd s}_{L_t^{\infty}L_x^2(I\times \R)} \lsm T N_0^{\frac 12+\alpha }\|V\|_{L_x^1}\|u\|_{L_t^{\infty}H_x^{\alpha}}^{2k+1}.
		}

	\end{lem}
\begin{proof}
		By Strichartz's estimates, the Sobolev and Kato-Ponce inequalities, for $\alpha>\frac 32$, we have
\begin{align*}
\normbb{\int_{0}^te^{i(t-s)\pd_x^2}\langle\partial_x\rangle^{\alpha} \big(\widetilde{V}G(u)\big)(s)\dd s}_{L_t^{\infty}L_x^2(I\times \R)} \lsm & \|\langle\partial_x\rangle^{\alpha}\big(\widetilde{V}G(u)\big)\|_{L_t^{1}L_x^{2}}\notag\\
\lsm &T\|\langle\partial_x\rangle^{\alpha} \widetilde{V}\|_{L_x^{\infty}}\|u\|_{L_{t}^{\infty}L_x^{2(2k+1)}}^{2k+1}\notag\\
&+T\| \widetilde{V}\|_{L_x^{\infty}}\|\langle\partial_x\rangle^{\alpha} u\|_{L_{t}^{\infty}L_x^2}\|u\|_{L_{t, x}^{\infty}}^{2k}\notag\\
\lsm &T\|\langle\partial_x\rangle^{\alpha} \widetilde{V}\|_{L_x^{\infty}}\|u\|_{L_{t}^{\infty}H_x^{\alpha}}^{2k+1}.
\end{align*}
This gives \eqref{easy-Vuu}.

For \eqref{easy-Vuu2}, it is easy to transfer the derivative to the high
frequency one. First, we have
$$
G(u)-G(u_{\ll N})=u_{\gtrsim N}\mathcal{O}(u, u_{\ll N}),
$$
where $\mathcal{O}(u, u_{\ll N})=\sum_{j=0}^{2k}u^{2k-j}u_{\ll N}^j$. Here we do not distinguish between $u$ and $\bar{u}$.

By the duality, Strichartz's estimates, local smoothing effect \eqref{Smooth2}, and Lemmas \ref{lem:Schur}, \ref{lem:littlewood-Paley}, we have
\begin{align*}
&\normbb{\sum_{N\in 2^{\N}}\int_0^t \langle\partial_x\rangle^{\alpha}e^{i(t-s)\pd_x^2}P_N\big(V(G(u)-G(u_{\ll N}))\big)(s)\dd s}_{L_t^{\infty}L_x^2(I\times \R)}\\
\lsm  &\sup_{h:\norm{h}_{L_x^2\leq 1}}\sum _{N\lsm M}\normbb{\Big\langle\int_{0}^t\langle\partial_x \rangle^{\alpha} e^{i(t-s)\pd_x^2}P_N(V P_Mu \mathcal{O}(u, u_{\ll N}) )\dd s, h\Big\rangle}_{L_t^{\infty}}\notag\\
\lsm &\sup_{h: \norm{h}_{L_x^2\leq 1}}\sum _{N\lsm M}\normbb{\int_{0}^t \langle\partial_x \rangle^{\alpha}  e^{i(t-s)\pd_x^2}P_N(V P_Mu \mathcal{O}(u, u_{\ll N})  )\dd s}_{L_t^{\infty}L_x^2}\|P_N h\|_{L_x^2}\notag\\
\lsm &\sup_{h: \norm{h}_{L_x^2\leq 1}}\sum _{N\lsm M}\frac{\langle N \rangle^{\alpha-\frac 12}}{\langle M \rangle^{\alpha -\frac 12}}\normbb{V \langle M \rangle^{\alpha-\frac 12 }P_Mu \mathcal{O}(u, u_{\ll N})  \dd s}_{L_x^{1}L_t^2}\|P_N h\|_{L_x^2}\notag\\
\lsm &\sup_{h: \norm{h}_{L_x^2\leq 1}}\sum _{N\lsm M}\frac{\langle N \rangle^{\alpha-\frac 12}}{\langle M \rangle^{\alpha-\frac 12 }}\|V\|_{L_x^1}\|\langle M \rangle^{\alpha-\frac 12 }P_Mu\|_{L_x^{\infty}L_t^2}\|\mathcal{O}(u, u_{\ll N} )\|_{L_{t, x}^{\infty}}\|P_N h\|_{L_x^2}\notag\\
\lsm &\|V\|_{L_x^1}\|\langle M \rangle^{\alpha-\frac 12 }P_Mu\|_{l_M^2L_x^{\infty}L_t^2}\|\mathcal{O}(u, u_{\ll N} )\|_{l_N^{\infty}L_{t, x}^{\infty}}\\
\lsm &T^{\frac 12}\|V\|_{L_x^1}\|\langle \partial_x \rangle^{\alpha }u\|_{L_t^{\infty}L_x^2}\|u\|_{L_t^{\infty}H_x^{\alpha}}^{2k}.
\end{align*}
Here the last step follows from Bernstein's inequality,
\begin{align*}
\|\langle M \rangle^{\alpha-\frac 12 }P_Mu\|_{l_M^2L_x^{\infty}L_t^2}\lsm \|\langle M \rangle^{\alpha}P_Mu\|_{l_M^2L_t^2L_x^{2}}
\lsm  T^{\frac 12}\|\langle \partial_x \rangle^{\alpha }u\|_{L_t^{\infty}L_x^2}.
\end{align*}
This gives \eqref{easy-Vuu2}. For \eqref{easy-Vuu3}, by Strichartz's estimates and Lemma \ref{lem:Bernstein},
\begin{align*}
&\normbb{\sum_{N: N< N_0}\int_0^t \langle\partial_x\rangle^{\alpha}e^{i(t-s)\pd_x^2}P_N\big(V|u_{\ll N}|^{2k}u_{\ll N}\big)\dd s}_{L_t^{\infty}L_x^2(I\times \R)}\\
\lsm & \sum_{N: N< N_0}N^{\alpha+\frac 12}\|V|u_{\ll N}|^{2k}u_{\ll N}\|_{L_{t, x}^1}\\
\lsm & T N_0^{\alpha +\frac 12}\|V\|_{L_x^1}\|u\|_{L_{t, x}^{\infty}}^{2k+1}\\
\lsm & T N_0^{\alpha+\frac 12 }\|V\|_{L_x^1}\|u\|_{L_t^{\infty}H_x^{\alpha}}^{2k+1}.
\end{align*}
This completes the proof of this lemma.
\end{proof}

	\subsection{Boundary terms}
	\begin{lem}\label{lem:nonlinear-estimate-boundary-1d-2}
		Let $I \subset\R^+$ be an interval containing $0$. Then, for any $N_0\in2^\N$,
		\EQn{\label{32511-12}
			\normb{e^{it\partial_x^2}\Om(\langle\partial_x\rangle^{\alpha-2}V, u_0)}_{L_t^{\infty}L_x^2(I\times \R)} \lsm \|P_{\geq N_0}\langle\partial_x\rangle^{\alpha-2}V\|_{L_x^{2}}\|u_0\|_{H_x^{\alpha}}^{2k+1},
		}
		and
		\EQn{\label{32512-12}
			\norm{\Om(\langle\partial_x\rangle^{\alpha-2}V, u)}_{L_t^{\infty}L_x^2(I\times \R)} \lsm  \|P_{\geq N_0}\langle\partial_x\rangle^{\alpha-2}V\|_{L_x^{2}}\|u\|_{L_t^{\infty}H_x^{\alpha}}^{2k+1}.
		}
	\end{lem}
	\begin{proof}
Recall the definition of $\Om$,
\EQ{
\Om(f_1,f)(x) := \sum_{N:N\geq N_0}\int_{\Ga(\vec\xi)} e^{ix\brkb{\sum_{j=1}^{2k+2}\xi_j}} m_N(\vec \xi) \wh{f_1}(\xi_1)\prod_{j=2}^{2k+2} \wh{f^{(\iota_j)}}(\xi_j) \dd \vec\xi.
}
Denote that
\EQn{\label{MN}
\Omega^N(f_1,f)(x) := \int_{\Ga(\vec\xi)} e^{ix\brkb{\sum_{j=1}^{2k+2}\xi_j}} m_N(\vec \xi) \wh{f_1}(\xi_1)\prod_{j=2}^{2k+2} \wh{f^{(\iota_j)}}(\xi_j) \dd \vec\xi.
}
Then
\EQn{\label{MN-1}
\Om(f_1,f)(x) := \sum_{N:N\geq N_0}\Omega^N(f_1,f)(x).
}
Since $m_N(\vec \xi)$ satisfies the conditions of Coifman-Meyer's multiplier, we have that, heuristically,
\EQ{
\Om^N(f_1,f)(x)\sim P_{\sim N} f_1\cdot  f_{\ll N}^{2k+1}.
}
Note that $\alpha>\frac 32$, by Strichartz's estimates, Lemma \ref{lem:Coifman-Meyer} and Sobolev's inequality, we obtain
		\begin{align}\label{32511-112}
			\normb{e^{it\partial_x^2}\Om(\langle\partial_x\rangle^{\alpha-2}V, u_0)}_{L_t^{\infty}L_x^2(I\times \R)}
\lsm &\normb{\Om(\langle\partial_x\rangle^{\alpha-2}V, u_0)}_{L_x^2}\notag\\
\lsm &\normb{P_{\sim N}\Om^N(\langle\partial_x\rangle^{\alpha-2}V, u_0)}_{L_x^2l_N^2(N\geq N_0)}\notag\\
\lsm & \normb{ \|P_N\langle\partial_x\rangle^{\alpha-2}V\|_{L_x^2}\|P_{\ll N}u_0\|_{L_x^{\infty}}^{2k+1}}_{l_N^2(N\geq N_0)}\notag\\
\lsm &  \|P_N\langle\partial_x\rangle^{\alpha-2}V\|_{L_x^2{l_N^2(N\geq N_0)}}\|P_{\ll N}u_0\|_{L_x^{\infty}l_N^{\infty}}^{2k+1}\notag\\
			\lsm &\|P_{\geq N_0} \langle\partial_x\rangle^{\alpha-2}V\|_{L_x^{2}}\|u\|_{L_t^{\infty}H_x^{\alpha}}^{2k+1},
		\end{align}
This gives \eqref{32511-12}. \eqref{32512-12} can be proved in the same way. Hence, we complete the proof of this lemma.
	\end{proof}

\subsection{High-order terms} Fix $0<\epsilon_0<\frac{1-b} {10}$. Since $b<1-\epsilon_0$, we may choose
$q_*$ so that
\begin{equation}\label{q-star-choice}
\frac1{1-\epsilon_0}<q_*<
\min\left\{\frac1b,\frac2{1-\epsilon_0}\right\}.
\end{equation}
Set $\theta_*:=(1-\epsilon_0)-\frac 1 {q_*}>0$. By Lemma
\ref{lem:potential-regularity}, $V\in L^{q_*}(\mathbb R)$. The upper bound on
$q_*$ is needed in the Bernstein estimate below. This restriction is
essential when $b<\frac 12$.
	\begin{lem}\label{lem:nonlinear-estimate-higher-order-1d-2}
		Let $I=[0, T) \subset\R^+$ with $T<1$. Then
		\EQ{
		&\sum_{l=2}^{2k+2}	\normbb{\int_0^t e^{i(t-s)\pd_x^2}\Om_l(\langle\partial_x\rangle^{\alpha-2}V, u, F(V, \widetilde{V}, u))(s)\dd s}_{L_t^{\infty}L_x^2(I\times \R)}\\
 \lsm& T^{\frac 12}\|\langle \partial_x \rangle^{\alpha-2} V\|_{L_x^{2}}\big(\|V\|_{L_x^{q_*}}+\|\widetilde{V}\|_{L_x^{\infty}}\big)\|u\|_{L_t^{\infty}H_x^{\alpha}}^{4k+1}.
		}
	\end{lem}
	\begin{proof}

By the smoothing effect \eqref{Smooth2}, we have
	\begin{align}\label{25-1201}
			&\normbb{\int_0^t e^{i(t-s)\pd_x^2}\Om_l(\langle\partial_x\rangle^{\alpha-2}V, u, F(V, \widetilde{V}, u))(s)\dd s}_{L_t^{\infty}L_x^2}\notag \\
			\lsm & \norm{\langle\partial_x \rangle^{-\frac 12}\Om_l(\langle\partial_x\rangle^{\alpha-2}V, u, F(V, \widetilde{V}, u))}_{L_x^1L_t^2}.
		\end{align}
		By Strichartz's estimates,
\begin{align}\label{25-1202}
			&\normbb{\int_0^t e^{i(t-s)\pd_x^2}\Om_l(\langle\partial_x\rangle^{\alpha-2}V, u, F(V, \widetilde{V}, u))(s)\dd s}_{L_t^{\infty}L_x^2}\notag \\
			\lsm & \norm{\Om_l(\langle\partial_x\rangle^{\alpha-2}V, u, F(V, \widetilde{V}, u))}_{L_{x, t}^{\frac 65}}.
		\end{align}
By \eqref{25-1201}, \eqref{25-1202}, and the interpolation, we have
\begin{align}\label{25-1203}
			&\normbb{\int_0^t e^{i(t-s)\pd_x^2}\Om_l(\langle\partial_x\rangle^{\alpha-2}V, u, F(V, \widetilde{V}, u))(s)\dd s}_{L_t^{\infty}L_x^2}\notag \\
			\lsm & \norm{\langle\partial_x \rangle^{-\frac 12(1-\epsilon_0)}\Om_l(\langle\partial_x\rangle^{\alpha-2}V, u, F(V, \widetilde{V}, u))}_{L_x^{r_0}L_t^{q_0}},
		\end{align}
where $r_0$ and $q_0$ satisfy $\frac 1{r_0}=1-\frac {\epsilon_0}{6}$ and $\frac 1 {q_0}=\frac 12+\frac {\epsilon_0}{3}$.

As in \eqref{MN-1}, we have
\EQn{\label{MN-2}
\Om_l(f_1,f,g)(x) := \sum_{N:N\geq N_0}\Omega_l^N(f_1,f, g)(x).
}
where $\Omega_l^N(f_1,f, g)(x)$ is defined as in \eqref {MN} with the function at the $l-th$ position replaced by $g$, where $2\leq l\leq 2k+2$.
Applying Lemma \ref{vector-estimate}, we have
\begin{align}\label{25-1204}
			&\norm{\langle\partial_x \rangle^{-\frac 12(1-\epsilon_0)}\Om_l(\langle\partial_x\rangle^{\alpha-2}V, u, F(V, \widetilde{V}, u))}_{L_x^{r_0}L_t^{q_0}}\notag \\
\lsm & \sum_{N:N\geq N_0} \norm{\langle\partial_x \rangle^{-\frac 12(1-\epsilon_0)}P_{\sim N}\Om_l^N(\langle\partial_x\rangle^{\alpha-2}V, u, F(V, \widetilde{V}, u))}_{L_x^{r_0}L_t^{q_0}}\notag \\
\lsm &  \sum_{N:N\geq N_0} N^{-\frac 12(1-\epsilon_0)}\|\langle\partial_x\rangle^{\alpha-2}P_N V\|_{L_x^2}\|u_{\ll N}\|_{L_{x, t}^{\frac {6k}{\epsilon_0}}}^{2k}\|P_{\ll N} F(V, \widetilde{V}, u)\|_{L_{x}^{\frac 2{1-\epsilon_0}}L_t^2}.
		\end{align}
Since $q_*<\frac 2{1-\epsilon_0}$, Lemma \ref{lem:Bernstein}, H\"older's and Sobolev's inequalities imply
\begin{align*}
\|P_{\ll N} \big(V|u|^{2k}u\big)\|_{L_{x}^{\frac 2{1-\epsilon_0}}L_t^2}\lsm & T^{\frac 12}\|P_{\ll N} \big(V|u|^{2k}u\big)\|_{L_{t}^{\infty}L_x^{\frac 2{1-\epsilon_0}}}\\
\lsm &T^{\frac 12} N^{\frac 1{q_*}-\frac{1-\epsilon_0}2}\|V|u|^{2k}u\|_{L_{t}^{\infty}L_x^{q_*}}\\
\lsm &T^{\frac 12} N^{\frac 1{q_*}-\frac{1-\epsilon_0}2}\|V\|_{L_x^{q_*}}\|u\|_{L_t^{\infty}H_x^{\alpha}}^{2k+1}.
\end{align*}
Moreover, we have
\begin{align*}
\|P_{\ll N} \big(\widetilde V|u|^{2k}u\big)\|_{L_{x}^{\frac 2{1-\epsilon_0}}L_t^2}\lsm & \|\widetilde{V}|u|^{2k}u\|_{L_{t}^2L_x^{\frac 2{1-\epsilon_0}}}\\
\lsm & T^{\frac 12} \|\widetilde{V}\|_{L_x^{\infty}}\|u\|_{L_t^{\infty}L_x^{2(2k+1)+}}^{2k+1}\\
\lsm & T^{\frac 12} \|\widetilde{V}\|_{L_x^{\infty}}\|u\|_{L_t^{\infty}H_x^{\alpha}}^{2k+1}.
\end{align*}
Combining the above three inequalities, using
$\|u_{\ll N}\|_{L_{x,t}^{\frac {6k}{\epsilon_0}}}
\lesssim\|u\|_{L_t^\infty H_x^\alpha}$ for $T<1$, we obtain
\begin{align}\label{25-1205}
			&\norm{\langle\partial_x \rangle^{-\frac 12(1-\epsilon_0)}\Om_l(\langle\partial_x\rangle^{\alpha-2}V, u, F(V, \widetilde{V}, u))}_{L_x^{r_0}L_t^{q_0}}\notag \\
\lesssim&T^{\frac12}\|\langle\partial_x\rangle^{\alpha-2}V\|_{L_x^2}
\left(N_0^{-\theta_*}\|V\|_{L_x^{q_*}}
+N_0^{-\frac12(1-\epsilon_0)}\|\widetilde V\|_{L_x^\infty}\right)
\|u\|_{L_t^\infty H_x^\alpha}^{4k+1}\notag \\
\lesssim&T^{\frac12}\|\langle\partial_x\rangle^{\alpha-2}V\|_{L_x^2}
\big(\|V\|_{L_x^{q_*}}+\|\widetilde V\|_{L_x^\infty}\big)
\|u\|_{L_t^\infty H_x^\alpha}^{4k+1},
		\end{align}
where $-\theta_*:=-(1-\epsilon_0)+\frac 1 {q_*}<0$.
Combining \eqref{25-1203} and \eqref{25-1205}, and then summing over the finitely many indices $l$, we finish the proof of this lemma.
\end{proof}

\subsection{Proof of Theorem \ref{thm-main}}
Based on the above several nonlinear estimates, we are in a position to prove Theorem \ref{thm-main}.
\begin{proof}[\bf{Proof of well-posedness in $H_x^{\alpha}$}]
Let $I=[0, T)\subset \R^+$ with $T<1$. First of all, by Strichartz's estimate, we have
		\EQn{\label{446}
			\norm{e^{it\pd_x^2} u_0}_{L_t^{\infty}H_x^{\alpha}(I\times \R)}=\|u_0\|_{H_x^{\alpha}}:=R.
		}
		Fix $0<\eta\ll 1$, by $V \in H_x^{\alpha-2}(\R)$, we take $N_0=N_0(V, \eta)\in 2^{\N}$, such that
		\EQn{\label{447}
			\|P_{\geq N_0}\langle\partial_x\rangle^{\alpha-2}V\|_{L_x^{2}}\leq \eta.
		}
		Denote the operator $\Phi$ by the following form,
		\EQ{
			\langle\partial_x\rangle^{\alpha}\Phi(u)=&\langle\partial_x\rangle^{\alpha}e^{i t\pd_x^2}u_0-i\lambda\int_{0}^t \langle\partial_x\rangle^{\alpha} e^{i(t-s)\pd_x^2}\widetilde{V}(x)|u|^{2k}u(s)\dd s\\
&-i\lambda\sum_{N\in 2^{\N}}\int_0^t \langle\partial_x\rangle^{\alpha}e^{i(t-s)\pd_x^2}P_N\big(V(|u|^{2k}u-|u_{\ll N}|^{2k}u_{\ll N})\big)\dd s\notag \\
    &-i\lambda\sum_{N: N< N_0}\int_0^t \langle\partial_x\rangle^{\alpha}e^{i(t-s)\pd_x^2}P_N\big(V|u_{\ll N}|^{2k}u_{\ll N}\big)\dd s\\
    &-\lambda \Om(\langle\partial_x\rangle^{\alpha-2}V, u)+\lambda e^{it\pd_x^2}\Om(\langle\partial_x\rangle^{\alpha-2}V, u_0)\notag\\
&+i\lambda^2 \sum_{l=2}^{2k+2}(-1)^{l-1}\int_0^t e^{i(t-s)\pd_x^2}\Om_l(\langle\partial_x\rangle^{\alpha-2}V, u, F(V, \widetilde{V}, u))(s)\dd s.
}
		Taking the working space as
		\EQ{
			B_{R}:=\{u\in C(I; H_x^{\alpha}(\R)):\|u\|_{L_t^{\infty}H_x^{\alpha}(I\times \R)}\leq 2R\}.
		}
We assume that $R\ge1$ without loss of generality.

		Next, we aim to prove $\Phi$ is a contraction in $B_{R}$. Hence, we need to collect the estimates of $\langle\partial_x\rangle^{\alpha}\Phi(u)$ in $L_t^{\infty}L_x^2$.
		
By Lemma \ref{lem:easy},
\EQn{\label{526-1}
		\normbb{\int_{0}^t\langle\partial_x\rangle^{\alpha} e^{i(t-s)\pd_x^2}\widetilde{V}(x)|u|^{2k}u(s)\dd s}_{L_t^{\infty}L_x^2(I\times \R)} \lsm T\|\langle\partial_x\rangle^{\alpha} \widetilde{V}\|_{L_x^{\infty}}R^{2k+1},
		}
\EQn{\label{526-2}
		\normbb{\sum_{N\in 2^{\N}}\int_0^t \langle\partial_x\rangle^{\alpha}e^{i(t-s)\pd_x^2}P_N\big(V(|u|^{2k}u-|u_{\ll N}|^{2k}u_{\ll N})\big)\dd s}_{L_t^{\infty}L_x^2(I\times \R)} \lsm T^{\frac 12}\|V\|_{L_x^1}R^{2k+1},
		}
and
\EQn{\label{526-3}
		\normbb{\sum_{N: N< N_0}\int_0^t \langle\partial_x\rangle^{\alpha}e^{i(t-s)\pd_x^2}P_N\big(V|u_{\ll N}|^{2k}u_{\ll N}\big)\dd s}_{L_t^{\infty}L_x^2(I\times \R)} \lsm T N_0^{\frac 12+\alpha }\|V\|_{L_x^1}R^{2k+1}.
		}
		By Lemma \ref{lem:nonlinear-estimate-boundary-1d-2},
			\EQn{\label{526-4}
			\normb{e^{it\partial_x^2}\Om(\langle\partial_x\rangle^{\alpha-2}V, u_0)}_{L_t^{\infty}L_x^2(I\times \R)} \lsm \|P_{\geq N_0}\langle\partial_x\rangle^{\alpha-2}V\|_{L_x^{2}}R^{2k+1},
		}
		and
		\EQn{\label{526-5}
			\norm{\Om(\langle\partial_x\rangle^{\alpha-2}V, u)}_{L_t^{\infty}L_x^2(I\times \R)} \lsm  \|P_{\geq N_0}\langle\partial_x\rangle^{\alpha-2}V\|_{L_x^{2}}R^{2k+1}.
		}
		By Lemma \ref{lem:nonlinear-estimate-higher-order-1d-2},
		\EQn{\label{526-6}
		&\sum_{l=2}^{2k+2}	\normbb{\int_0^t e^{i(t-s)\pd_x^2}\Om_l(\langle\partial_x\rangle^{\alpha-2}V, u, F(V, \widetilde{V}, u))(s)\dd s}_{L_t^{\infty}L_x^2(I\times \R)}\\
 \lsm& T^{\frac 12}\|\langle \partial_x \rangle^{\alpha-2} V\|_{L_x^{2}}\big(\|V\|_{L_x^{q_*}}+\|\widetilde{V}\|_{L_x^{\infty}}\big)R^{4k+1}.
		}
		By the estimates \eqref{526-1}-\eqref{526-6}, for $T<1$, there exists a constant $C$ depending on $\|V\|_{H_x^{\alpha-2}}$, $\|V\|_{L_x^1}$, $\|V\|_{L_x^{q_*}}$, $\|\langle \partial_x \rangle^{\alpha}\widetilde{V}\|_{L_x^{\infty}}$, such that
		\begin{align}\label{448}
			\norm{\langle\partial_x\rangle^{\alpha }\Phi(u)}_{L_t^{\infty}L_x^2(I\times \R)}\leq  R+C\eta R^{2k+1}+CT^{\frac 12}R^{4k+1}+CT N_0^{\frac 12+\alpha}R^{2k+1}.
		\end{align}
		First, by \eqref{447}, we take $N_0=N_0(\eta, \|V\|_{H_x^{\alpha-2}})$ large to obtain a small $\eta$, such that
		\EQ{
			C\eta R^{2k}\leq \frac 14.
		}
	Then, we take $T=T(R, N_0, \|V\|_{H_x^{\alpha-2}}, \|V\|_{L_x^1}, \|V\|_{L_x^{q_*}}, \|\langle \partial_x \rangle^{\alpha}\widetilde{V}\|_{L_x^{\infty}})$ small enough so that
\EQ{
CT^{\frac 12}R^{4k}+C T N_0^{\frac 12+\alpha}R^{2k}\leq \frac 14.
}
		Therefore, by the above estimates, we have
		\EQ{
			\norm{\Phi(u)}_{L_t^{\infty}H_x^{\alpha }(I\times \R)}\leq 2R.
		}
		Hence, we have that $\Phi:B_{R}\rightarrow B_{R}$. The contraction mapping then follows using a similar and standard argument. Therefore, this gives the existence and uniqueness of $u\in C([0, T);H_x^{\alpha}(\R))$ to the equation \eqref{eq:NLS}, proving part \textit{(1)} of Theorem \ref{thm-main}.
\end{proof}

	\begin{proof}[\bf{Proof of ill-posedness}]
		Next, we prove that for $\alpha \geq \frac 52-b$, the equation \eqref{eq:NLS} is ill-posed in $H_x^{\alpha}(\R)$. We only need to show the ill-posedness in $H_x^{\frac 52-b}$. This shall be achieved by an iteration scheme based on the Duhamel formula, see Bejenaru and Tao \cite{B-T-06}.
If the data-to-solution map is $C^{2k+1}$ at the origin as a map into $C([0,T];H^{\frac52-b}(\R))$, then its $(2k+1)$-st derivative at the origin yields the necessary estimate
$$
\Big\|\int_0^t e^{i(t-s)\partial_x^2}(|x|^{-b} |e^{is\partial_x^2}u_0|^{2k}e^{is\partial_x^2}u_0)\dd s\Big\|_{H_x^{\frac 52-b}(\R)}\lsm \|u_0\|_{H_x^{\frac 52-b}(\R)}^{2k+1}.
$$
For convenience, we denote
\EQ{
			A [u_0]=\int_0^t e^{-is\partial_x^2}(|x|^{-b} |e^{is\partial_x^2}u_0|^{2k}e^{is\partial_x^2}u_0)\dd s.
		}
To achieve our goal, it is sufficient to show that for any $T>0$ and $M_0>0$, there exists $u_0\in H_x^{\frac 52-b}$ such that
\begin{align}\label{1211-1}
\sup\limits_{t\in [0,T]}\|A[u_0](t)\|_{H_x^{\frac 52-b}}\geq M_0.
\end{align}
We choose the initial data $u_0$ such that
		$$u_0=P_{\leq 1}u_0, \quad \widehat{u}_0\geq 0, \mbox{ and } \int_{\R} \widehat{u}_0(\xi)\dd \xi>0.$$
		(For example, $u_0(x):=P_{\leq 1}e^{-|x|^2}$ satisfies the above conditions).
	
Moreover, let $N_0$ be a large constant determined later, and
		\EQ{
			t:=\frac 1{N_0},
		}
		and the set
		\EQ{
			\Omega:=\bigcup_{l=N_0}^{N_0^2}\Omega_l:=\bigcup_{l=N_0}^{N_0^2}\Big(\sqrt{N_0}\sqrt{2l\pi+\frac{5\pi}{12}}, \sqrt{N_0}\sqrt{2l\pi+\frac{\pi}{2}}\Big).
		}
		Notice that $\Omega_l\cap \Omega_j=\varnothing$ if $l\neq j$.
We also have that
\begin{align*}
\|A [u_0]\|_{H_x^{\frac52-b}(\R)}= \norm{\langle\xi\rangle^{\frac 52-b }\widehat{A[u_0]}(\xi)}_{L_{\xi}^2(\R)}\geq \norm{\langle\xi\rangle^{\frac 52-b }\widehat{A[u_0]}(\xi)}_{L_{\xi}^2(\Omega)}.
\end{align*}
Next, we restrict our analysis to the frequency region $\Omega$. This choice corresponds to the situation where all components of the initial data are localized at low frequencies,
 while the nonlinear interaction generates a high-frequency output concentrated in $\Omega$. Roughly speaking, the nonlinear interaction generates frequencies of size at least
 $N_0$, i.e. $|\xi|\geq N_0$.

For $A[u_0]$, using the Fourier transform of the distribution $|x|^{-b}$ and the choice of $u_0$, we have
		\EQ{
			\widehat{A [u_0]}(\xi)=&\,c_b\int_0^t\int_{\xi=\sum_{j=1}^{2k+2}\xi_j}e^{is \Phi(\vec \xi)}|\xi_1|^{b-1}\prod_{j=2}^{2k+2}\widehat{u_0^{(\iota_j)}}(\xi_j)\dd \xi_1\cdots \dd \xi_{2k+1}\dd s\\
=&\,c_b\int_0^t\int_{\xi=\sum_{j=1}^{2k+2}\xi_j}e^{is \Phi(\vec \xi)}|\xi_1|^{b-1}\prod_{j=2}^{2k+2}\widehat{u_0}^{(\iota_j)}(\iota_j\xi_j)\dd \xi_1\cdots \dd \xi_{2k+1} \dd s,
		}
where $\iota_j$ is defined in Definition \ref{defn:normal-form}, and \EQ{
\Phi(\vec \xi):=|\xi|^2- \sum_{j=2}^{2k+2}\iota_j \xi_j^2.
}
		For $|\xi|\geq N_0$ with $N_0$ large enough, and note that $|\xi_j|\leq 2$ ($j=2, 3, \ldots, 2k+2$) by the definition of $u_0$, the phase function $\Phi(\vec \xi)$ does not vanish.
Applying the integration-by-parts, we have
		\EQ{
			\widehat{A [u_0]}(\xi)
			=&\,c_b\int_{\xi=\sum_{j=1}^{2k+2}\xi_j}\frac{e^{it \Phi(\vec \xi)}-1}{i\Phi(\vec \xi)}|\xi_1|^{b-1}\prod_{j=2}^{2k+2}\widehat{u_0}^{(\iota_j)}(\iota_j\xi_j)\dd \xi_1\cdots \dd \xi_{2k+1}.
		}
Taking the real part of $\widehat{A [u_0]}(\xi)$ yields
\EQn{\label{reA1010}
		\re	\big(\widehat{A [u_0]}(\xi)\big)
			=\,c_b\int_{\xi=\sum_{j=1}^{2k+2}\xi_j}\frac{\sin[t\Phi(\vec \xi)]}{\Phi(\vec \xi)}|\xi_1|^{b-1}\prod_{j=2}^{2k+2}\widehat{u_0}^{(\iota_j)}(\iota_j\xi_j)\dd \xi_1\cdots \dd \xi_{2k+1}.
		}
By the mean value theorem, we have
\EQ{
\sin[t\Phi(\vec \xi)]=\sin (t|\xi|^2)+O(t\cdot \sum_{j=2}^{2k+2}\iota_j \xi_j^2).
}
where the notation $O(f)$ means that there exists a constant $0<C\leq 1$, independent of $t$ and $\vec \xi$,
such that $|O(f)|\leq C|f|\leq |f|$.

Noting that if $\xi \in \Omega$, then  $t|\xi|^2\in (2l\pi+\frac{5\pi}{12}, 2l\pi+\frac{\pi}{2})$, which further implies $\sin (t|\xi|^2)\ge \frac 12$. Moreover, for $N_0$ large enough, we have $|t\cdot \sum_{j=2}^{2k+2}\iota_j \xi_j^2|\leq \frac {4(2k+1)}{N_0} \leq \frac 1 4$. Hence, we conclude that
\begin{align}\label{1010}
\sin[t\Phi(\vec \xi)]\geq \frac 14.
\end{align}	
By the estimates \eqref{reA1010} and \eqref{1010}, we obtain
\EQn{\label{reA1010-1}
		\re	\big(\widehat{A [u_0]}(\xi)\big)
			\geq &\,C\int_{\xi=\sum_{j=1}^{2k+2}\xi_j}\frac{|\xi_1|^{b-1}}{\Phi(\vec \xi)}\prod_{j=2}^{2k+2}\widehat{u_0}^{(\iota_j)}(\iota_j\xi_j)\dd \xi_1\cdots \dd \xi_{2k+1}\\
\geq &\,C\int_{\xi=\sum_{j=1}^{2k+2}\xi_j}\frac{1}{|\xi|^{3-b}}\prod_{j=2}^{2k+2}\widehat{u_0}^{(\iota_j)}(\iota_j\xi_j)\dd \xi_1\cdots \dd \xi_{2k+1}.
		}	
	Noting that $\mbox{Re}\big(\widehat{A [u_0]}(\xi)\big)>0$, the above inequality yields that
\begin{align*}
\|A [u_0]\|_{H_x^{\frac52-b}(\R)}\geq \norm{\langle\xi\rangle^{\frac 52-b }\widehat{A[u_0]}(\xi)}_{L_{\xi}^2(\Omega)}\geq \norm{\langle\xi\rangle^{\frac 52-b }\mbox{Re}\big(\widehat{A[u_0]}(\xi)\big)}_{L_{\xi}^2(\Omega)}.
\end{align*}	
Hence, we conclude that
\begin{align}\label{es-A}
\|A [u_0]\|_{H_x^{\frac 52-b}(\R)}^2\geq \int_{\Omega}\big|\langle\xi\rangle^{\frac 52-b}\mbox{Re}\big(\widehat{A[u_0]}(\xi)\big)\big|^2\dd \xi
\geq  C_0\int_{\Omega}\frac {1}{|\xi|}\dd \xi,
		\end{align}
		where $C_0:=C(\int \widehat{u}_0(\xi)d\xi)>0$.
		
Recall that $\Omega_l\cap \Omega_j=\varnothing$ if $l\neq j$, we have
		\EQ{
			\int_{\Omega}\frac {1}{|\xi|}\dd \xi=&\sum_{l=N_0}^{N_0^2}\int_{\sqrt{N_0}\sqrt{2l\pi+\frac{5\pi}{12}}}^{\sqrt{N_0}\sqrt{2l\pi+\frac{\pi}{2}}}\frac {1}{|\xi|}\dd \xi\\
			=&\frac 12\sum_{l=N_0}^{N_0^2}\ln\frac {2l\pi+\frac{\pi}{2}}{2l\pi+\frac{5\pi}{12}}.
		}
		Notice that $\ln(1+x)\geq \frac 12x$ for $0\leq x\leq  1$. Hence, taking $N_0$ large enough, for any $l\geq N_0$,
		\EQ{
			\ln\frac {2l\pi+\frac{\pi}{2}}{2l\pi+\frac{5\pi}{12}}\gtrsim \frac {1}{2l\pi+\frac{5\pi}{12}}.
		}
		Further, by the above two estimates, we have
		\begin{align}\label{6311}
			\int_{\Omega}\frac {1}{|\xi|}\dd \xi\gtrsim &  \sum_{l=N_0}^{N_0^2} \frac {1}{2l\pi+\frac{5\pi}{12}}\notag\\
			\gtrsim & \ln \frac {2N_0^2\pi+\frac{5\pi}{12}}{2N_0\pi+\frac{5\pi}{12}}\notag\\
			\gtrsim & \,\ln N_0.
		\end{align}
		Hence, by the estimates \eqref{es-A} and \eqref{6311}, we obtain \eqref{1211-1}.
		Therefore, the proof of ill-posedness is complete. This proves part \textit{(2)} and finishes the proof of Theorem \ref{thm-main}.
	\end{proof}

	\vskip 0.2cm

\end{document}